\documentclass [a4paper,english]{article}
\usepackage{amsmath, amsfonts, amssymb}
\usepackage[cp1251]{inputenc}
\usepackage[T2A]{fontenc}
\usepackage{babel,indentfirst,epigraph}
\usepackage{graphicx,graphics,verbatim}
\newcommand{\bn}[1]{\[#1\]}
\newcommand{\be}[2]{\begin{equation}\label{#1} {#2}\end{equation}}

\title{A conjecture on monotonicity of a ratio of Kummer hypergeometric functions}
\author{Sitnik\,S.M.}

\begin{document}
\maketitle

\begin{abstract}
In the preprint  \cite{Sit1} the author formulated some conjectures on  monotonicity of  ratios for exponential series remainders. They are equivalent to  conjectures on monotonicity of a ratio of Kummer hypergeometric functions and presumably not proved still. In this short note the most interesting conjecture from   \cite{Sit1} is reproduced with its generalizations.
\end{abstract}

Let us consider exponential series remainder in the form
\be{rem}{R_n(x)=\exp(x)-\sum_{k=0}^{k=n}\frac{x^k}{k!}=\sum_{k=n+1}^{k=\infty}\frac{x^k}{k!}, \,x\ge 0.}
Different inequalities for $R_n(x)$ and its combinations were proved by H.\,Hardy, W.\,Gautschi, H.\,Alzer and many others \cite{Sit1}.

In the preprint  \cite{Sit1} the author studied inequalities of the form
\be{inq}{m(n)\le f_n(x)=\frac{R_{n-1}(x)R_{n+1}(x)}{\left[R_n(x)\right]^2} \le M(n), \,x\ge 0.}
The search for best constants $m(n)=m_{best}(n),\,M(n)=M_{best}(n)$ has some history. The left--hand side of (\ref{inq}) was first proved by Kesava Menon in 1943 with
$m(n)=\frac{1}{2}$ (not best) and by Horst Alzer in 1990 with $m_{best}(n)=\frac{n+1}{n+2}=f_n(0)$ \cite{Sit1}. In fact this result is a special case of an inequality proved by Walter Gautschi in 1982 in \cite{Gau1} (cf. \cite{Sit1}).

It seems that the right--hand side of (\ref{inq}) was first proved by the author in \cite{Sit1} with $M_{best}=1=f_n(\infty)$. In \cite{Sit1} dozens of generalizations of inequality   (\ref{inq}) and related results were proved.

Obviously the above inequalities are consequences of the next conjecture  formulated in \cite{Sit1}.

\textit{\textbf{Conjecture 1.} The function $f_n(x)$ in  (\ref{inq}) is monotone increasing for $x\in [0; \infty), n\in \mathbb{N}$.}

In 1990's we tried to prove this conjecture by considering $(f_n(x))^{'}\ge 0$ and expanding  triple products of  hypergeometric functions but failed \cite{Sit2}.

Consider a representation via Kummer hypergeometric functions
\be{Kum}{f_n(x)=\frac{n+1}{n+2} \ g_n(x), g_n(x)=\frac{_1F_1(1; n+1; x)  _1F_1(1; n+3; x)}{\left[_1F_1(1; n+2; x)\right]^2}.}
So the conjecture 1 is equivalent to  the next conjecture 2.

\textit{\textbf{Conjecture 2.} The function $g_n(x)$ in  (\ref{Kum}) is monotone increasing for $x\in [0; \infty), n\in \mathbb{N}$.}

This leads us to the next more general

\textit{\textbf{Problem 1.} Find monotonicity  in $x$ conditions for $x\in [0; \infty)$ for all parameters a,b,c for the function}
\be{prob}{h(a,b,c,x)=\frac{_1F_1(a; b-c; x)  _1F_1(a; b+c; x)}{\left[_1F_1(a; b; x)\right]^2}.}

We may also call (\ref{prob}) \textit{The abc--problem}, why not?

Another generalization is to change Kummer  hypergeometric functions to higher ones.

\textit{\textbf{Problem 2.} Find monotonicity  in $x$ conditions for $x\in [0; \infty)$ for all vector--valued parameters a,b,c for the function}
\be{probpq}{h_{p,q}(a,b,c,x)=\frac{_pF_q(a; b-c; x)  _pF_q(a; b+c; x)}{\left[_pF_q(a; b; x)\right]^2},}
\bn{a=(a_1,\ldots, a_p), b=(b_1,\ldots, b_q), c=(c_1, \ldots, c_q).}

Similar problems were recently carefully studied by D.\,Karp and his coauthors, cf. \cite{Karp1}--\cite{Karp10}.
These problems are also connected with famous Tur\'{a}n--type inequalities \cite{Bar1}--\cite{Bar4}. Cf. also \cite{SS}.

\bigskip

\textsl{\textbf{
I am glad to announce that the above conjectures 1, 2 and problems 1, 2 were proved in 2014 by the author and Khaled Mehrez (Tunisia).}} 

\bigskip

The proofs are published in \cite{MS1}--\cite{MS2}. All above mentioned conjectures and problems were also generalized to $q$--hypergeometric functions in  \cite{MS3}--\cite{MS4}.

\bigskip

Sergei M. Sitnik.\\ Voronezh Institute of the Russian Ministry of Internal Affairs.
Voronezh, Russia.\\
email: pochtaname@gmail.com

\end{document}